\documentclass[12pt,a4paper]{amsart}

\usepackage[a4paper]{geometry}
\usepackage{amsmath,amsthm,amsfonts,amssymb,mathrsfs,amsbsy}

\usepackage{hyperref}
\usepackage{mathtools}
\mathtoolsset{showonlyrefs}

\newtheorem{theorem}{Theorem}

\newtheorem{definition}{Definition}
\newtheorem{lemma}{Lemma}

\title[Asymptotic Behavior of a Fractional Volterra Equations]{Asymptotic behavior of solutions of a time-space fractional diffusive Volterra equation}

\author[S. Ahmad]{Sofwah Ahmad} 
\address[S. Ahmad]{Department of Mathematics, College of Computing and Mathematical Sciences, Khalifa University of Science and Technology, P.O. Box 127788, Abu Dhabi, UAE\\ 
\href{https://orcid.org/0000-0001-7641-7759}{orcid.org/0000-0001-7641-7759}}
\email{100059797@ku.ac.ae}

\author[M. Kirane]{Mokhtar Kirane} 
\address[M. Kirane]{Department of Mathematics, College of Computing and Mathematical Sciences, Khalifa University of Science and Technology, P.O. Box 127788, Abu Dhabi, UAE\\ 
\href{https://orcid.org/0000-0002-4867-7542}{orcid.org/0000-0002-4867-7542}}
\email{mokhtar.kirane@ku.ac.ae}
\urladdr {https://www.ku.ac.ae/college-people/mokhtar-kirane}

\begin{document}

\begin{abstract}
    In this paper, we study the time-space fractional differential equation of the Volterra type:
    \begin{align*}
        \mathcal{D}^\alpha_{0 \vert t} (u) +(-\Delta_N)^{\sigma}u &= u(1+au-bu^2)-au\int_0^t\mathcal{K}(t-s) u(\cdot) \, ds,
    \end{align*}
    where $a,b>0$ are given constants, $\alpha,\sigma \in (0,1)$, equipped with a homogeneous Neumann's boundary condition and a positive initial data.  The boundedness and uniform continuity of the solution on the entire $\mathbb{R}^+$ are established.
    Moreover, the asymptotic behavior of the positive solution is investigated.

\medskip
    
    
    \medskip
    
    \noindent \textbf{MSC Classification (2020):} 26A33, 35B09, 35B40, 34K37, 35R09, 35R11.

\end{abstract}

\maketitle

    \section{\textbf{Introduction}} \label{sec1}
	In this paper, we are interested to study the time asymptotic behavior of the bounded solution of the following problem 
		\begin{align}
			\mathcal{D}^\alpha_{0 \vert t} (u) +(-\Delta_N)^{\sigma}u &= u(1+au-bu^2)-au\mathcal{K}* u,  \label{Eq1} \\
			u(x,0) &= u_0(x), \label{Eq2}\;   
		\end{align}
	posed on $Q=\Omega\times(0,T)$ where $\Omega\subset\mathbb{R}^d$ (where $d=1,2$ in practice) is a bounded domain with smooth boundary $\partial\Omega$ with outward normal $\eta$; the asterisk ($*$) denotes the time convolution, i.e
		\begin{equation*}
			\mathcal{K}* u(x,t) = \int_{0}^{t} \mathcal{K}(t-s)u(\cdot,s) ds;
		\end{equation*}
	the constants $a$ and $b$ in the problem are assumed to be positive and $\mathcal{K}$ is the delay kernel function satisfying
		\begin{equation*}
			{\mathcal{K}} \ge 0,\  
                {\mathcal{K}} \in C^1(0,\infty)\cap L^1(0,\infty),\ 
                \int_{0}^{\infty} {\mathcal{K}}(t)dt = 1.
		\end{equation*}
	The derivative operator $\mathcal{D}^\alpha_{0 \vert t}$ is the Caputo fractional derivative of order $0\le \alpha \le 1$, the operator $(-\Delta_N)^\sigma, 0<\sigma<1$, is the fractional power of the classical Laplacian operator with its specified domain and $u_0(x)\ge 0$ is a given bounded function.
	
	Equation \eqref{Eq1} describes the evolution of a single species population whose density is $u$ at time $t$ in the bounded domain $\Omega$. As in the existing model with classical derivatives, the first term in the right-hand side represents the non-delay term with growth rate $a$ due to advantages of aggregating as a group (e.g protected against predator and increased reproduction success) and the decay rate $b$ is due to inter-species competition; the delay represented by convolution term in the Volterra equation \eqref{Eq1} represents the regulatory term with the assumption that population growth depends on a smooth average over the past history of the population rather than only information at a certain instantaneous time. More discussion about various population models can be found in \cite{Britton1} and \cite{Cushing}. 
	
	In this paper, we use the time fractional derivative $\mathcal{D}^\alpha_{0 \vert t} u$ rather than the usual derivative to represent the accumulation of the information on the population density variation within the specified time interval $(0,T)$; the fractional Laplacian term $(-\Delta_N)^{\sigma}u$ describes spatial distribution of the species with its possibility to move to remote sites inside the region $\Omega$ within a certain period of time (see \cite{pozrikidis2016fractional}). With the fractional Laplacian, it is implied that the  model under consideration is equipped with the Neumann's boundary condition 
		\begin{equation}
			\nabla u\cdot \eta = 0,  \;  (x,t) \in \partial \Omega\times(0,T), \label{Eq3}\\
		\end{equation}
	where $\eta$ is the outer normal of $\partial\Omega$, which says that no species migration into or out of $\Omega$.
	
	Before proceeding further, let us dwell on the existing literature on the classical differential equations related to our problem (i.e  $\alpha=\sigma=1$). In \cite{Britton2}, the classical differential equation version of \eqref{Eq1} is studied with the delay term $(1+a-b)u\mathcal{K}* u$ instead of $au\mathcal{K}* u$, where the authors showed that for $0<b<1+a$ and $4a<(3(b-1)+\sqrt{b^2+6b+1})$, $u(x,t)$ goes to 1 uniformly for $x\in\bar{\Omega}$. Schiaffino \cite{Schiaffino}, Redlinger \cite{Redlinger2}, and Yamada \cite{yamada1}, studied the same model with intraspecific competion term $u\mathcal{K} * u$ subject to Neumann's boundary condition.  The asymptotic properties of solutions of the problem with Dirichlet's boundary condition is also studied by Yamada \cite{yamada2}. In Yoshida and Yamada \cite{yamada3}, the study of the spatial inhomogeneity is taken into consideration. Kirane, Guedda, and Tatar \cite{GueddaKiraneTatar1997} also studied problem \eqref{Eq1}; by means of an entropy functional, they showed that the solution tends to $1/\sqrt{b}$ as $t\to\infty$.

    In this work, using the method of Langlais-Phillips, that is outlined in the paper of Kirane \cite{Kirane1989}, by noting that the orbits of the dynamical system are precompact in $C^\mu(\Omega)$, for some $\mu>0$ that will be precised later, we show that the solution of the problem tends to $1/\sqrt{b}$ which is the positive solution of the limiting problem of \eqref{Eq1}-\eqref{Eq2}.

\section{\textbf{Preliminaries}}
	\subsection{Definition and notation}

	First, we give the definitions of Riemann-Liouville fractional integral,  Riemann-Liouville fractional derivative and Caputo fractional derivative (see e.g. \cite{KilbasSrivastavaTrujillo}) 
        \begin{definition}
            Let $f\in L^{1}([0,T])$ and $\alpha>0$.  The left and right Riemann-Liouville fractional integrals of order $\alpha>0$ are defined by
                \begin{align}
                        (I_{0|t}^\alpha f)(t) &:= \frac{1}{\Gamma(\alpha)} \int_0^t \frac{f(\tau) }{(t - \tau)^{1-\alpha}} \, d\tau, 
                        \\
                        (I_{t|T}^\alpha f)(t) &:= \frac{1}{\Gamma(\alpha)} \int_t^T \frac{f(\tau)}{(\tau - t)^{1-\alpha}} \, d\tau,
                \end{align}
            respectively.
        \end{definition}

		\begin{definition}
			The Riemann-Liouville fractional derivative  of order $ \alpha \in (0,1]$ of an absolutely continuous function $f$ is defined by
                \begin{align}
                    ({}^{RL}D_{0|t}^\alpha y)(t) &:= \frac{d}{dt} (I_{0|t}^{1-\alpha} y)(t) 
                    = \frac{1}{\Gamma(1 - \alpha)} \frac{d}{dt} \int_0^t \frac{y(\tau) }{(t - \tau)^{\alpha}}\, d\tau, 
                    \quad (t > 0),
                    \\[10pt]
                    ({}^{RL}D_{t|T}^\alpha y)(t) &:= -\frac{d}{dt} (I_{t|T}^{1-\alpha} y)(t) 
                    = \frac{-1}{\Gamma(1 - \alpha)} \frac{d}{dt} \int_t^T \frac{y(\tau)}{(\tau - t)^{\alpha}}\, d\tau, 
                    \quad (t < T).
                \end{align}

		\end{definition}

         The following relation between Riemann-Liouville integral and derivative of order $0 < \alpha \leq 1$ holds:
         \begin{lemma}\cite[Lemma 2.4]{KilbasSrivastavaTrujillo}
            If $\alpha > 0$ and $f(t) \in L^p(0, T)$ ($1 \le p \le \infty$), then the following equalities
                \begin{align}
                    ({}^{RL}D_{0|t}^\alpha I_{0|t}^\alpha f)(t) &= f(t), \quad \text{and} \quad ({}^{RL}D_{t|T}^\alpha I_{t|T}^\alpha f)(t) = f(t), 
                \end{align}
            hold almost everywhere on $[0, T]$.
        \end{lemma}
        
		\begin{definition}
			The Caputo fractional derivative  of order $0 < \alpha \leq 1$ of a differentiable function $f$ is defined by
				\begin{equation}
					({\mathcal{D}}_{0 \vert t}^{\alpha} f)(t) := \frac{1}{\Gamma(1-\alpha)} \int_{0}^{t} (t-\tau)^{- \alpha} f'(\tau) d\tau 
                    \quad 0\le t\le T.
				\end{equation}
		\end{definition}

	Next, we consider the definition of the fractional Laplacian given via the spectrum of the classical Laplacian operator $-\Delta$ on the same domain \cite{gal2020fractional}. 
	Let $\{\lambda_n\}_{n=0}^{\infty}$ with $\lambda_0=0<\lambda_1<\lambda_2<\ldots$ be the sequence of eigenvalues of the operator $-\Delta$ on $\Omega$ that is equipped with Neumann boundary condition on $\partial\Omega$,  and $\{e_n\}_{n=0}^\infty$ be the corresponding normalized eigenfunctions in $\Omega$, i.e, 
	\begin{eqnarray*}
		-\Delta e_n &= \lambda_n e_n, &x\in\Omega,\\
            \frac{\partial e_n}{\partial \eta} &= 0, \quad &x\in\partial\Omega.
	\end{eqnarray*}
    Consider the space 
        \begin{equation}
		H^\sigma(\Omega):= \left\{u=\sum_{n=0}^{\infty} u_n e_n \in L^2(\Omega); 
            \quad \|u\|_{H^\sigma(\Omega)}=\sum_{n=1}^{\infty}\lambda_n^\sigma u_n^2 < +\infty \right\}.
        \end{equation}
    Defining the inner product
        \begin{equation*}
            \langle u,v \rangle_{H^\sigma(\Omega)} := \sum_{n=1}^{\infty}\lambda_n^\sigma u_n v_n , 
            \qquad\text{ for } \qquad 
            u=\sum_{n=0}^{\infty} u_ne_n,\ v= \sum_{n=0}^{\infty} v_ne_n,
	\end{equation*}
    we see that the space $H^\sigma(\Omega)$, equipped with the inner product $\langle u,v \rangle_{H^\sigma(\Omega)}$, is a Hilbert space.

        \begin{definition}
		The fractional Laplacian operator $(-\Delta_N)^\sigma$ is defined by 
                \begin{equation}
                     (-\Delta_N)^\sigma u(x) := \sum_{n=1}^{\infty} \lambda_n^\sigma u_n e_n,
                \end{equation}
            with domain
                \begin{equation}
                   D \big((-\Delta_N)^\sigma\big) = \left\{u\in H^\sigma(\Omega); \quad \frac{\partial u}{\partial \eta}=0 \right\}.
                \end{equation}
    	\end{definition}
	
    \begin{definition}
        The \textit{two parameters Mittag-Leffler function} $E_{\alpha,\beta}(z)$ is given by
            \begin{equation}
                E_{\alpha,\beta}(z) = \sum_{k=0}^\infty \frac{z^k}{\Gamma(\alpha k +\beta)} \quad (z,\beta\in \mathbb{C};\ Re(\alpha)>0).
            \end{equation}
    \end{definition}

	\subsection{Preliminary lemmas}
	In this subsection, we give the preliminary lemmas that will be used in proving our result.
            \begin{lemma} \cite[Lemma 2.22]{KilbasSrivastavaTrujillo}
			For a continuous function $f(t)$, the following relationship between fractional integral and Caputo fractional derivative holds 
			\begin{equation}
					I_{0|t}^\alpha \mathcal{D}_{0|t}^{\alpha } f(t) = f(t) - f(0) \label{eqn J^aD^a}.
                \end{equation}
		\end{lemma}

		\begin{lemma}\label{Lemma_chain_convexD^a}
			If $\varphi \in C^1(\mathbb{R})$ is a convex function and $x\in C([0,T])\cap C^1((0,T])$, for some $T>0$, then 
				\begin{equation*}
					\mathcal{D}^\alpha_{0 \vert t} \varphi(x)(t) \le \varphi'(x(t))\mathcal{D}^{\alpha}_{0 \vert t} x(t),
					\qquad0<\alpha\le1.
				\end{equation*}
		\end{lemma}
            \begin{proof}
                Since $\varphi$ is convex function, we have 
                    \begin{equation*}
                        \varphi(x)-\varphi(y) - \varphi'(y)(x-y)\ge \frac{1}{2}\varphi''(\lambda x + \mu y )(x-y)^2, 
                        \quad\text{where }\quad 0\le\lambda,\mu<1.
                    \end{equation*}
                Define
                    \begin{equation*}
                        \Psi(s,t)= \varphi(x(s))-\varphi(x(t))-\varphi'(x(t))\left(x(s)-x(t)\right),
                    \end{equation*}
                then    
                    \begin{equation*}
                        \frac{\partial\Psi(s,t)}{\partial s} =  \varphi'(x(s)) x'(s)-\varphi'(x(t))x'(s).
                    \end{equation*}

                Thus 
                \begin{align*}
                    \mathcal{D}^{\alpha}_{0 \vert t} \varphi(x)(t) &- \varphi'(x(t))\mathcal{D}^{\alpha}_{0 \vert t} x(t)
                    \\
                    &=\frac{1}{\Gamma(1-\alpha)}\left(\int_0^t \frac{\varphi'(x(s))\ x'(s)}{(t-s)^\alpha}ds-\varphi'(x(t)) \int_0^t \frac{x'(s)}{(t-s)^\alpha}ds\right)\\
                    &= \frac{1}{\Gamma(1-\alpha)}\left(\int_0^t \frac{\left(\varphi'(x(s))-\varphi'(x(t))\right)x'(s)}{(t-s)^\alpha} ds \right)\\
                    &=\frac{1}{\Gamma(1-\alpha)} \left(\int_0^t \frac{\partial\Psi(s,t)}{\partial s}(t-s)^{-\alpha}  ds\right)\\
                    &= \frac{1}{\Gamma(1-\alpha)} \left(\lim_{s\to t} \frac{\Psi(s,t)}{(t-s)^\alpha}-\frac{\Psi(0,t)}{t^\alpha}-\alpha \int_0^t \frac{\Psi(s,t)}{(t-s)^{\alpha+1}}ds\right)\\
                    &= \frac{1}{\Gamma(1-\alpha)} \left(\lim_{s\to t} \frac{\partial\Psi(s,t)/\partial s}{\alpha (t-s)^{\alpha-1}} -\frac{\Psi(0,t)}{t^\alpha}-\alpha \int_0^t \frac{\Psi(s,t)}{(t-s)^{\alpha+1}}ds\right)\\
                    &= \frac{1}{\Gamma(1-\alpha)} \left( -\frac{\Psi(0,t)}{t^\alpha}-\alpha \int_0^t \frac{\Psi(s,t)}{(t-s)^{\alpha+1}}ds\right),
                \end{align*}
                and the result follows.
            \end{proof}
  
	 \textbf{Remark.} In the case $\varphi(u)=u^p$, we obtain the generalisation of J.I. Diaz, T. Pierantozi and L. Vazquez conjectured inequality, which was obtained by Alsaedi, Ahmad and Kirane \cite{AAK2}.

     We need the following result on weakly singular Gronwall's inequality.     
        \begin{lemma} \cite[Lemma 6.19]{diethelm2010analysis}
        \label{ineq-weakly Gronwall}
            Let $\alpha,C,T\in \mathbb{R}^+$. If $\psi:[0,T]\to \mathbb{R}$ is a non-negative continuous function that satisfies
                \begin{equation*}
                    \psi(t) \le \psi(0)+\frac{C}{\Gamma(\alpha)} \int_0^t (t-s)^{\alpha-1} \psi(s) ds,
                \end{equation*}
            for all $t\in [0,T]$, then
                \begin{equation*}
                    \psi(t) \le \psi(0) E_{\alpha,1}(Ct^\alpha),
                \end{equation*}
            where $E_{\alpha,\beta}(z)$ is the two parameters Mittag-Leffler function. 
        \end{lemma}

    Using the properties of scalar products, the following integration by parts formula for the fractional Laplacian immediately follows.
        \begin{lemma}
            For $0<\sigma\le 1$ and $u,v\in D\left((-\Delta_N)^\sigma\right)$, it holds 
                \begin{align}
                    \int_{\Omega} u (-\Delta_N)^\sigma v \ dx =  \int_{\Omega} v (-\Delta_N)^\sigma u \ dx.
                \end{align}
            \label{Lemma_int by parts lap}
        \end{lemma}
    We also need the following inequality of Stroock and Varopoulos.
        \begin{lemma}\cite[Theorem 1]{LiskevichSemenov}
            Let $0<\sigma\le 2$ and $p> 1$. For nonnegative $u\in L^p(\Omega)$ with $(-\Delta_N)^{\frac{\sigma}{2}}u \in L^p(\Omega)$, we have
                \begin{align}
                    \int_{\Omega} u^{p-1} (-\Delta_N)^{\frac{\sigma}{2}}u \ dx 
                    \ge \frac{4(p-1)}{p^2} \int_{\Omega} \left|(-\Delta_N)^{\frac{\sigma}{4}} u^{\frac{p}{2}}\right|^2 dx. 
                \end{align}
            \label{Lemma_Strook&Varopulous_ineq}
        \end{lemma}

    We also have the following relation for time-fractional integration by parts. 
    \begin{lemma}\cite[Theorem 11]{Almeida_2019}
       Let $ 0 < \alpha < 1 $. The following relations hold:
            \begin{align}
                \int_0^T y(t) \, \mathcal{D}_{0|t}^{\alpha} x(t) \, dt 
                &= \int_0^T x(t) \,  {}^{RL}D_{t|T}^{\alpha} y(t) \, dt + \left[ x(t) \, I_{t|T}^{1 - \alpha} y(t) \right]_{t = 0}^{t = T}. 
                \label{eq: frac_int_by_parts1}
            \end{align} 
    \label{lemma: int by parts}
    \end{lemma}

\section{\textbf{Local Existence}}

 In this section, we prove the existence and uniqueness of local (in time) solution $u(x,t)$ on some interval $[0,T]$, for $T<\infty$. Let $$f(x,t,u) = u(1+au-bu^2)-au\mathcal{K}* u.$$ Furthermore, let the operators $S_\alpha(t), P_\alpha(t):H^\sigma(\Omega)\to H^\sigma(\Omega)$ be defined as
     \begin{align} 
        & S_\alpha(t) v:=\int_0^{\infty} \Phi_\alpha(\tau) S\left(\tau t^\alpha\right) v \ d \tau, \label{eqnSa2}\\
        \intertext{and}
        & P_\alpha(t) v:=\alpha t^{\alpha-1} \int_0^{\infty} \tau \Phi_\alpha(\tau) S\left(\tau t^\alpha\right) v \ d \tau, \label{eqnPa2}
    \end{align}
with $S(t)$ is  the strongly continuous semigroup generated by the operator $A:=(-\Delta_N)^\sigma$ and $\Phi_\alpha(z)$ is the Wright type function (see Kilbas et al. \cite{KilbasSrivastavaTrujillo}) 
    \begin{equation*}
        \Phi_\alpha(z):=\sum_{n=0}^{\infty} \frac{(-z)^n}{n ! \Gamma(-\alpha n+1-\alpha)}, \quad 0<\alpha<1, \quad z \in \mathbb{C}.
    \end{equation*}
Then, we have the following definition for the representation of the solution of the problem \eqref{Eq1}-\eqref{Eq3}.

\begin{definition}
    Let $u_0\in L^\infty(\Omega)$ and $T>0$. We say that $u\in C([0,T],L^\infty(\Omega))$ is a mild solution of the problem \eqref{Eq1}-\eqref{Eq3} if it can be expressed as 
        \begin{equation}\label{eq-solution}
            u(t) = S_\alpha(t)u_0 + \int_0^t P_\alpha (t-s)f(s) \ ds,
        \end{equation}
    for $t\in [0,T]$.
\end{definition}

It is useful to first notice the boundedness of operators $S_{\alpha}$ and $P_{\alpha}$.

\begin{lemma}\cite[Proposition 2.1.3]{gal2020fractional}
    There exists constants $C_s,C_p>0$ such that for all $t>0$ and $v\in H^{\sigma}(\Omega)$, we have
        \begin{enumerate}
            \item $\|S_{\alpha}(t)v\|_{H^{\sigma}(\Omega)} \le C_s \|v\|_{H^{\sigma}(\Omega)}$,
            \item[] 
            \item $\|t^{1-\alpha}P_{\alpha}(t)v\|_{H^{\sigma}(\Omega)} \le C_p \|v\|_{H^{\sigma}(\Omega)}$.
        \end{enumerate}
    \label{Lemma_S&P-bound}
\end{lemma}

\begin{theorem}(Local existence)
     Let $u_0\in C(\Omega)$; then there exists a maximal time $T_{max}>0$ and a unique mild solution $u\in C([0, T_{max}); L^\infty(\Omega))$ of the problem \eqref{Eq1}-\eqref{Eq2} with the alternative
        \begin{itemize}
            \item either $T_{max}=+\infty$;
            \item or $T_{max}<+\infty$ and
            $$\lim_{t\to T_{max}}\|u(t)\|_\infty=\lim_{t\to T_{max}} \sup_{0<t<T_{max}}|u(t)|=+\infty. $$
        \end{itemize}
 \end{theorem}
    \begin{proof}
    We prove the local existence through Banach's fixed point theorem. Let $u_0:\Omega \to \mathbb{R}^+$ be such that $\|u_0\|_{L^\infty(\Omega)} \equiv \rho_0$ small enough. 
    Choose $\rho$ such that $\rho - C_s \rho_0>0$ and let 
        \begin{equation*}
            B_\rho:=\{u\in C([0,T];L^\infty(\Omega)): \|u\|_{B_\rho}:=\displaystyle\sup_{t\in [0,T]} \|u(\cdot,t)\|_{L^\infty(\Omega)} \le\rho \}.
        \end{equation*}
    
    We define the operator $F$ on $B_\rho$ as 
        \begin{equation}
            (Fu)(t) :=  S_\alpha (t) u_0 + \int_0^t P_{\alpha} (t-s) (f)(x,s,u(s)) \ ds.
        \end{equation}
    To use Banach's fixed point theorem, we shall prove that $F$ is a bounded operator and it is a contraction on $B_\rho$.

    \begin{itemize}
        \item We first show that $F:B_\rho \to B_\rho$. With the aid of Lemma \ref{Lemma_S&P-bound}, we have
            \begin{align*}
                \|Fu\|_{L^\infty(\Omega)} &\le \|S_\alpha(t) u_0\|_{L^\infty(\Omega)} + \int_0^t \|P_{\alpha} (t-s)\| \|f(u(s))\|_{L^\infty(\Omega)} \ ds\\
                &\le C_s \|u_0\|_{L^\infty(\Omega)} + \int_0^t (t-s)^{\alpha -1} C_p \|f(u)\|_{L^\infty(\Omega)}\ ds.
            \end{align*}
        We note that
            \begin{align*}
                \|f(u)\|_{L^\infty(\Omega)} &\le \|u +  au^2-bu^3\|_{L^\infty(\Omega)}+ \bigg\|au \int_0^t K(t-s) u(x,s) ds\bigg\|_{L^\infty(\Omega)} \\
                &\le \|u\|_{B_\rho} + a\|u\|_{B_\rho}^2+b\|u\|_{B_\rho}^3 + a \|u\|_{B_\rho}^2 \int_0^\infty \mathcal{K}(s) ds\\
                &\le \rho + 2a\rho^2 + b\rho^3,
            \end{align*}
        so 
            \begin{align*}
                \|Fu\|_{B_\rho}  &\le  C_s\rho_0 + C_p (\rho+2a\rho^2+b\rho^3) \displaystyle\sup_{t\in [0,T]} \left\{\int_0^t (t-s)^{\alpha-1} \ ds \right\}\\
                &\le C_s\rho_0 + C_p (\rho +2a\rho^2+b\rho^3) \displaystyle\sup_{t\in [0,T]} \left\{\frac{t^\alpha}{\alpha}\right\}\\
                &\le C_s\rho_0 + C_p (\rho+2a\rho^2+b\rho^3) \frac{T^\alpha}{\alpha}\\
                &\le \rho,
            \end{align*}
        by choosing $T\le T_1$ where 
            \begin{equation}
                T_1 = \left[ \frac{\alpha(\rho-C_s\rho_0)}{C_p(\rho+2a\rho+b\rho^2)} \right]^{1/\alpha}.
            \end{equation}
        Therefore, we conclude that $F:B_\rho\to B_\rho$. 
        
        \item Next we show that $F$ is a contraction on $B_\rho$.
        
        We observe that $f$ is Lipschitz continuous on $B_\rho$ since for any $u,v\in B_\rho$, we have 
            \begin{align*}
                f(u)-f(v) &= (u-v) + a(u^2-v^2) -  b(u^3-v^3)\\
                &\quad -a\left[(u-v)(\mathcal{K}*u)-v(\mathcal{K}*(u-v))\right]\\
                &= (u-v) (1+ a(u+v) -  b(u^2+uv+v^2)) \\
                &\quad-a\left[(u-v)(\mathcal{K}*u)-v(\mathcal{K}*(u-v))\right],
            \end{align*}
        so
            \begin{align*}
                \|f(u)-f(v)\|_{B_\rho}  
                &\le \|u-v\|_{B_\rho} \|1+ a(u+v) -  b(u^2+uv+v^2)\|_{B_\rho} \\
                & \quad + a \|u-v\|_{B_\rho} \|\mathcal{K}*u\|_{B_\rho}+ a\|v\|_{B_\rho}\|\mathcal{K}*(u-v)\|_{B_\rho}.
            \end{align*}
        Whereupon,
            \begin{align*}
                \|f(u)-f(v)\|_{B_\rho} 
                &\le \|u-v\|_{B_\rho}  \left(1+ a(2\rho) +  b(3\rho^2) +2 a\rho  \int_0^\infty \mathcal{K}(s)\ ds \right)\\
                &:=C_f\|u-v\|_{B_\rho}. 
            \end{align*}
        Therefore,
        \begin{align*}
                \|Fu-Fv\|_{B_\rho}  &= \bigg \| \int_0^t  P_{\alpha}(t-s) (f(u(s)-f(v(s))\ ds\bigg \|_{B_\rho} \\
                &\le C_p\ C_f \|u-v\|_{B_\rho}  \displaystyle\sup_{t\in [0,T]} \left\{\int_0^t (t-s)^{\alpha-1} \ ds \right\}\\
                &\le C_p \ C_f \|u-v\|_{B_\rho}  \frac{T^\alpha}{\alpha}\\
                &\le \frac{1}{2} \|u-v\|_{B_\rho},
            \end{align*}    
        for $T\le T_2$ where $T_2=\left[\frac{\alpha}{2C_p \ C_f} \right]^{1/\alpha}$. By letting $T_3=\min\{T_1,T_2\}$, we see that $F$ is a contraction for $t\in[0, T_3]$ and hence, by Banach's fixed point theorem, the problem has a unique solution on $[0,T_3]$.
    
        \item Let
            \begin{align*}
                T_{\max} := \sup \{T>0: \text{ there exists a mild solution } u\in B_{\rho}\}\le+\infty.
            \end{align*}
        We show that either $T_{\max}=+\infty$, or $T_{\max}<+\infty$ and (in the latter case)  $\|u\|_{L^\infty(\Omega)} \to +\infty$ as $t\to T_{\max}$. If the former holds, we say that the solution is global. To prove the latter part, we proceed by contradiction. Suppose that $u$ is a solution on the interval $[0,T]$ with $T<T_{\max}=+\infty$ but $\|u\|_{L^\infty(\Omega)}<\infty$. For some $\tau>0$, we define $w\in C((0, T_{\max}+\tau),C(\Omega))$ as 
            \begin{align*}
            w(t)=
                \begin{cases}
                    u(t), &t\in[0,T_{\max}],\\
                    u(t-\tau), &t\in[T_{\max}, T_{\max}+\tau].
                \end{cases}
            \end{align*}
        We see that $w(t)$ is also a mild solution of the problem \eqref{Eq1}-\eqref{Eq2}. This contradicts the definition of $T_{\max}$.  
        
        \item We prove that the solution is unique on the whole interval of existence. 
        
        Let $u,\bar{u}\in B_\rho$ be two mild solutions of  problem \eqref{Eq1}-\eqref{Eq3} on $[0,T_{\max})$, then
            \begin{align*}
                \|u-\bar{u}\|_{L^\infty(\Omega)}\le \int_0^t  \|P_{\alpha}(t-s) (f(u(s)-f(\bar{u}(s)) \|_{L^\infty(\Omega)}\ ds
                \\
                \le C_p\ C_f \int_0^t (t-s)^{\alpha-1}\ \|u(s)-\bar{u}(s) \|_{L^\infty(\Omega)}\ ds.
            \end{align*}
        Let $\Psi(t)=\|u(\cdot,t)-\bar{u}(\cdot,t)\|_{L^\infty(\Omega)}$. Note that $\Psi(0)=0$. By Lemma \ref{ineq-weakly Gronwall}, we conclude that 
        $$u-\bar{u}=0 \quad\Longleftrightarrow\quad u=\bar{u}.$$ 
    \end{itemize}
\end{proof}

\section{\textbf{Positivity, Boundedness and Uniform Continuity of the Solution} }
We have established that the solution is bounded on  
$0<t<T<T_{\max}$. In this section, we show that the solution is positive and bounded for all $t\in\mathbb{R}$.	
	\begin{theorem}\label{positive&bounded}
		Assuming $u_0\ge 0,u_0\not\equiv 0$ and $u_0\in C(\bar{\Omega})$, the solution of problem \eqref{Eq1}-\eqref{Eq2} satisfies 
			\begin{equation*}
				0\le u(x,t)\le \max\{\|u_0\|_{\infty}, R\} \text{ on }\bar{\Omega}\times(0,\infty),
			\end{equation*}
		where $R$ is the positive solution of $1+aw-bw^2=0$. Observe that $R\ge \frac{1}{\sqrt{b}}$.
	\end{theorem}
    
    \begin{proof}
        First we prove the positivity. Let $u^+=max(u,0)$ and $u^- = min(u,0)$. Write $u=u^+ + u^-$. We note that $u^+u^-=0$, $u^2 = (u^+ + u^-)(u^+ + u^-)=(u^+)^2 + (u^-)^2$, and $u^3 = u\cdot u^2= (u^+)^3 + (u^-)^3$.
    
        Taking the scalar product of \eqref{Eq1} by $u^-$ in $L^2(\Omega)$, we obtain
            \begin{align*}
               \int_{\Omega} u^- \mathcal{D}^\alpha_{0 \vert t}  u^- + u^-(-\Delta_N)^{\sigma}(u^-) \ dx
               &= \int_{\Omega}(u^-)^2 + a(u^-)^3-b(u^-)^4\ dx
               \\
               &\quad -a\int_{\Omega}(u^-)^2\mathcal{K}*u^- \ dx.
               \end{align*}
        Using integration by parts, we get
            \begin{align*}
                \int_{\Omega} u^- \mathcal{D}^\alpha_{0 \vert t}  u^- \ dx
                &=  -\int_\Omega |(-\Delta_N)^{(\sigma/2)}(u^-)|^2\ dx 
                \\
                &\quad + \int_{\Omega}(u^-)^2 + a(u^-)^3-b(u^-)^4
               \\
               &\quad -a\int_{\Omega}(u^-)^2\mathcal{K}*u^- \ dx;
            \end{align*}
        which implies
            \begin{align*}
               \int_{\Omega} u^- \mathcal{D}^\alpha_{0 \vert t}  u^- dx
               \le \int_{\Omega}(u^-)^2 -a(u^-)^2\mathcal{K}*u^- \ dx.
            \end{align*}
        Since the local existence is established, there exists $T_{max}$ such that $u$ satisfies the fractional differential equation \eqref{Eq1} on $[0,T_{max})$. This means that the solution $u$ is bounded on $[0,T]$ for all $T<T_{max}$, i.e. there exists a continuous function $C(t)$ such that $\|u(t)\|_{L^\infty(\Omega)}\le C(t)$. Thus
            \begin{align*}
               \int_{\Omega} u^- \mathcal{D}^\alpha_{0 \vert t}  u^- \ dx
               &\le \int_{\Omega}(u^-)^2 +aC (u^-)^2  \left(\int_0^t K(s)ds\right) \ dx
               \\
               &\le \int_{\Omega}(u^-)^2 +aC (u^-)^2 \left(\int_0^\infty K(s)ds\right) \ dx
               \\
               &\le \big(1+a C\big)\int_{\Omega}(u^-)^2 \ dx.
            \end{align*}
        where $C= \displaystyle\sup_{t\in [0,T]} C(t)$. Since by Lemma \ref{Lemma_chain_convexD^a}, $\frac{1}{2} \mathcal{D}^\alpha_{0 \vert t}  (u^-)^2 \le u^- \mathcal{D}^\alpha_{0 \vert t}  u^-$, we then have
            \begin{align*}
               \frac{1}{2} \int_{\Omega} \mathcal{D}^\alpha_{0 \vert t}  (u^-)^2 dx 
               \le  \big(1+aC\big)\int_{\Omega}(u^-)^2dx.  
            \end{align*}
      
        Let $\psi(t)=\int_{\Omega}(u^-)^2$. As $u(x,0)=u_0\ge 0$, then $\psi(0)=\int_\Omega u^-(x,0)^2 dx=0$. By exchanging the integral and $\mathcal{D}_{0|t}^\alpha$ on left-hand side, we get 
            \begin{align}\label{eqn psi}
                \begin{cases}
                    \mathcal{D}^\alpha_{0 \vert t} \psi(t) &\le 2(1+aC) \psi(t),\\
                    \psi(0) &= 0.
                \end{cases}
            \end{align}
        By applying \eqref{eqn J^aD^a}, we have 
            \begin{align*}
                \psi(t) &\le \psi(0)+ 2(1+aC) I_{0|t}^\alpha \ \psi(t) ,
            \end{align*}
        which by Lemma \ref{ineq-weakly Gronwall} leads to
            \begin{equation*}
                \psi(t) = 0     
                \quad\Longleftrightarrow\quad 
                \int_\Omega (u^-)^2 dx = 0
                \quad\Longrightarrow\quad 
                u^- = 0
                \quad\Longrightarrow\quad 
                u=u^+\ge 0.
            \end{equation*}
            
        To prove the boundedness of $u$, we consider all possible cases. If $0\le u\le R$ then the proof is done. Otherwise let $u\ge R$. Multiplying equation \eqref{Eq1} by $u^{p-1}$ and integrating over $\Omega$, we get 
            \begin{align}
                \int_\Omega u^{p-1} \mathcal{D}^\alpha_{0 \vert t} u\ dx 
                &+\int_\Omega u^{p-1} (-\Delta_N)^{\sigma}u \ dx
                =\int_\Omega u^{p-1} f(u)\ dx.
                \label{eq: boundedness f(u)}
                \intertext{For $u>R$, we have $u+au^2-bu^3<0$. As for the convolution term we have }
                \mathcal{K}*u\ge \mathcal{K}*R 
                \,\,\Longrightarrow\,\, 
                 u\mathcal{K}*u
                 &\ge R\mathcal{K}*R 
                 \,\,\Longrightarrow\,\, 
                 -au\mathcal{K}*u\le -aR\mathcal{K}*R <0.
                \intertext{Hence, by using Lemma \ref{Lemma_chain_convexD^a} and \eqref{Lemma_Strook&Varopulous_ineq}, we obtain}
                \frac{1}{p} \int_\Omega \ \mathcal{D}^\alpha_{0 \vert t} (u^{p})\ dx 
                &+\frac{4(p-1)}{p^2} \int_{\Omega} \left|(-\Delta_N)^{\frac{\sigma}{2}} u^{\frac{p}{2}}\right|^2  \ dx
                \le 0;
                \intertext{whereupon,}
                \int_\Omega \ \mathcal{D}^\alpha_{0 \vert t} (u^{p})\ dx &\le 0.
                \intertext{Integrating with respect to $t$ $\alpha$-times and using \eqref{eqn J^aD^a}, we obtain}
                \int_\Omega \  u^{p}\ dx 
                &\le \int_\Omega \ u_0^{p}\ dx.
            \end{align}
        Taking $p\to \infty$, we get $\|u\|_{L^\infty(\Omega)}\le \|u_0\|_{L^\infty(\Omega)}$.
    \end{proof}
    
    \subsection{Uniform Continuity of $u$}
    To prove the uniform continuity of the solution of \eqref{Eq1}-\eqref{Eq3},  we notice that $A$ is a sectorial operator (see Gal and Warma \cite{gal2020fractional}). The following lemma shall be useful (see Guswanto and Suzuki \cite{Guswanto}).
    
    \begin{lemma}\cite[Theorems 2.1, 2.2]{Guswanto}
        Let $A$ be a sectorial operator and let $S_\alpha(t)$ and $P_\alpha(t)$ be the operators defined by \eqref{eqnSa2} and \eqref{eqnPa2}, respectively. Then, the following statements hold:
            \begin{enumerate}
                \item[(i)]  $S_\alpha(t), P_\alpha(t) \in B((0,\infty) ; H^\sigma(\Omega))$ for $t>0$, and if $x \in H^\sigma(\Omega)$ then 
                    \begin{equation}
                        A S_\alpha(t) x = S_\alpha(t) A x, \quad
                        A P_\alpha(t) x = P_\alpha(t) A x;\label{eqn_AP=PA}
                    \end{equation}
                \item[(ii)] the functions $t\mapsto S_\alpha(t)$, $t\mapsto P_\alpha(t)$ belong to $C^\infty((0,\infty);B(H^\sigma(\Omega)))$ ;
                \item[(iii)]for $x\in H^\sigma(\Omega)$ and $t>0$,
                    \begin{equation}
                        \frac{d}{dt} S_\alpha(t) x = AP_\alpha(t) x.
                    \label{eqn_S'=AP}
                    \end{equation}
            \end{enumerate}
    \end{lemma}
    
    For $0<\beta<1$, we consider the fractional power $A^{-\beta}$ (see Henry \cite{Henry})
        \begin{equation*}
            A^{-\beta} = \frac{1}{\Gamma(\beta)}\int_0^t t^{\beta-1} e^{-At} dt,
        \end{equation*}
    and $A^\beta$ is considered as the inverse of $A^{-\beta}$ with $D(A^\beta) = R(A^{-\beta})$. Moreover, for a fixed $p>d$, the following embedding holds true:
        \begin{equation}\label{embedding_1}
            D(A^\beta)\hookrightarrow C^\nu (\Omega), \qquad \ \text{if}\ \ 0\le \nu < 2\beta - \frac{d}{p};
        \end{equation}
    with $2\beta - \frac{d}{p}<1$.

    In the next lemma, we consider some estimates involving the fractional  operator $A^\beta$ and the families of operators  $\left\{S_\alpha(t)\right\}_{t>0},\left\{P_\alpha(t)\right\}_{t>0}$ (see Guswanto and Suzuki \cite{Guswanto}). 
        \begin{lemma}\cite[Theorem 2.6]{Guswanto}
            For each $0<\beta<1$, there exist positive constants $C_1=C_1(\alpha, \beta)$ and $C_2=C_2(\alpha, \beta)$ such that, for all $x \in H^\sigma(\Omega)$,
            \begin{align}
                \left\|A^\beta S_\alpha(t) x\right\|_{H^\sigma(\Omega)} 
                &\leq C_1 t^{-\alpha}\left(t^{\alpha(1-\beta)}+1\right)\|x\|_{H^\sigma(\Omega)}, \quad\ t>0, \label{eqn_boundSa}
                \\
                \left\|A^\beta P_\alpha(t) x\right\|_{H^\sigma(\Omega)} 
                &\leq C_2 t^{\alpha(1-\beta)-1}\|x\|_{H^\sigma(\Omega)}, 
                \qquad\qquad\quad t>0 . \label{eqn_boundPa}
            \end{align}
        \end{lemma}

    \begin{theorem}
        For any $\delta >0$, the function $t\mapsto u(x,t)$ is uniformly continuous with respect to $t\in[\delta,\infty)$ in the topology of $L^2(\Omega)$.
        \label{uniform_cont thm}
    \end{theorem}
    \begin{proof}
        Applying $A^\beta$ to $u(x,t+h)-u(x,t)$ of \eqref{eq-solution} for any $h>0$ gives
            \begin{align*}
                A^\beta u(x,t+h)-A^\beta u(x,t) 
                &= A^\beta\left(S_{\alpha}(t+h)-S_{\alpha}(t)\right) u_0 
                \\
                &\quad + \int_0^{t+h} A^\beta P_{\alpha} (t+h-\tau) f(x,\tau,u(\tau)) \ d\tau
                \\
                &\quad - \int_0^{t} A^\beta P_{\alpha} (t-\tau) f(x,\tau,u(\tau)) \ d\tau
                \\
                &= A^\beta\left(S_{\alpha}(t+h)-S_{\alpha}(t)\right) u_0 
                \\
                &\quad + \int_t^{t+h} A^\beta P_{\alpha} (t+h-\tau) f(x,\tau,u(\tau)) \ d\tau 
                \\
                & \quad - \int_0^{t} A^\beta \left(P_{\alpha} (t+h-\tau)-P_{\alpha} (t-\tau)\right) f(x,\tau,u(\tau)) \ d\tau.
            \end{align*}
        The first term tends to 0 as $h\to 0$, since 
        \begin{align*}
            \|A^\beta\left(S_{\alpha}(t+h)-S_{\alpha}(t)\right) u_0\|_{L^2(\Omega)}
            &= \|A^\beta h S'_\alpha \big(t+\delta h\big) u_0 \|_{L^2(\Omega)}\\
            &= h \|A^\beta A P_\alpha \big(t+\delta h\big) u_0\|_{L^2(\Omega)}\\
            &= h \|A^\beta P_\alpha \big(t+\delta h\big) A u_0\|_{L^2(\Omega)}\\
            &= h C_2(\alpha,\beta) \big(t+\delta h\big) \|A u_0\|_{L^2(\Omega)},
        \end{align*}
    for $\delta\in(0,1)$, where we have used \eqref{eqn_AP=PA}, \eqref{eqn_S'=AP} and \eqref{eqn_boundPa}, respectively.

    For the second term, we observe that
        \begin{align*}
            \bigg\|\int_t^{t+h} A^\beta P_\alpha(t+h-\tau)\ &f(u(\tau)) \ d\tau \bigg\|_{L^2(\Omega)}
            \\
            &= \bigg\|\int_0^{h} A^\beta P_\alpha(\tau) f(u(t+h-\tau)) \ d\tau \bigg\|_{L^2(\Omega)}
            \\
            & \le \int_0^{h}\big\| A^\beta P_\alpha(\tau) f(u(t+h-\tau)) \big\|_{L^2(\Omega)} \ d\tau. 
        \end{align*}
    Note that, by Theorem \ref{positive&bounded}, $u$ is bounded by $R$. Hence, for $t>0$, we have 
        \begin{align*}
            \|f(u(\tau)\|_{L^2(\Omega)} 
            &\le \|u +  au^2-bu^3\|_{L^2(\Omega)} + \bigg\|au \int_0^t K(t-s) u(x,s) ds\bigg\|_{L^2(\Omega)} 
            \\
            &\le R + 2aR^2 + bR^3:=R^*.
        \end{align*}
    Therefore, using \eqref{eqn_boundPa}, we have 
        \begin{align*}
            &\int_0^{h}\big\| A^\beta P_\alpha(\tau) f(u(t+h-\tau)) \big\|_{L^2(\Omega)} \ d\tau \\
            &\quad\le \int_0^h C_2\tau^{\alpha(1-\beta)-1} \big\| f(u(t+h-\tau)) \big\|_{L^2(\Omega)} d\tau
            \\
            &\quad \le \frac{C_2 R^*}{\alpha(1-\beta)}t^{\alpha(1-\beta)}\bigg|_0^h
            =\frac{C_2 R^*}{\alpha(1-\beta)}h^{\alpha(1-\beta)},
        \end{align*}
    which tends to 0 as $h\to 0$.

    Since 
        \begin{align*}
            &\big\|A^\beta \left(P_{\alpha} (t+h-\tau)-P_{\alpha} (t-\tau)\right)f(u(\tau))\big\|_{L^2(\Omega)}
            \\
            &\quad\le \big\|A^\beta P_{\alpha} (t+h-\tau) f(u(\tau))\big\|+\big\|A^\beta P_{\alpha}(t-\tau) f(u(\tau))\big\|_{L^2(\Omega)} 
            \\
            &\quad\le C_2R^* \left((t+h-\tau)^{-\alpha(\beta-1)-1}+(t-\tau)^{-\alpha(\beta-1)-1}\right),
        \end{align*}
    and for any $h>0$, we have 
        \begin{align*}
            C_2R^* &\int_0^t (t+h-\tau)^{-\alpha(\beta-1)-1}+(t-\tau)^{-\alpha(\beta-1)-1} \ d\tau
            \\
            &=\frac{C_2 R^*}{\alpha(1-\beta)} \left[(t+h-\tau)^{\alpha(1-\beta)} + (t-\tau)^{\alpha(1-\beta)}\right]_0^t\\
            &\le\frac{C_2 R^*}{\alpha(1-\beta)} \left[ h^{\alpha(1-\beta)}+(t+h)^{\alpha(1-\beta)} + t^{\alpha(1-\beta)} \right] \\
            &< \infty.
        \end{align*}
    Thus, by the dominated convergence theorem, we have 
        \begin{equation}\label{eqn-LDCT}
            \lim_{h\to 0} \big\|\int_0^t A^\beta \left(P_{\alpha} (t+h-\tau)-P_{\alpha} (t-\tau)\right)f(u(\tau)) \ d\tau\big\|_{L^2(\Omega)} = 0.
        \end{equation}
    In view of \eqref{embedding_1} and \eqref{eqn-LDCT}, we have
        \begin{equation}
            \lim_{h\to 0}\|u(x,t+h)-u(x,t)\|_{L^2(\Omega)}=0.
        \end{equation}
     \end{proof}

\section{\textbf{Asymptotic stability of the positive equilibrium}}
 
In this section, we prove the global attractivity of the positive equilibrium. We denote by $\omega(u_0)$ the $\omega$-limit set for $u_0$:
    \begin{equation}
        \omega(u_0) = \{\varphi\in H^\sigma(\Omega)\cap L^{\infty}(\Omega) 
        \ s.t\
        \exists t_n \to \infty, \ u(\cdot,t_n)\to \varphi(\cdot) \text{ in } L^2(\Omega)\}.
    \end{equation}
    
    Due to the compact embedding $C^\nu(\Omega) \hookrightarrow C^\mu(\Omega)$ for $0<\mu<\nu<1$, and the boundedness of the solution, the set $\omega(u_0)$ is nonempty as the orbit $\displaystyle\bigcup_{ t\ge 0}\{u(t)\}$ is precompact in $C^{\mu}(\bar\Omega)$. We shall see that any non-negative $\varphi\in\omega(u_0)$ solves the limiting problem associated to problem \eqref{Eq1}-\eqref{Eq2}, namely
    \begin{align}
        (-\Delta_N)^\sigma \varphi &= \varphi(1-b\varphi^2) \quad\text{in }  \Omega.
        \label{eq: stationary}
    \end{align}

Our result is contained in the following theorem:

\begin{theorem}
    {
    Let $u\in L^{\infty}(\Omega\times (0,1))$ be a weak solution of the problem \eqref{Eq1}-\eqref{Eq2}. Then $u$ tends to the constant stationary solution of the elliptic problem \eqref{eq: stationary} uniformly, that is, for any $\varphi\in \omega(u_0)$
        \begin{align*}
            (-\Delta_N)^\sigma \varphi &= \varphi(1-b\varphi^2) \quad\text{in }  \Omega.
        \end{align*}
    Moreover, $\varphi = 1/\sqrt{b}$.
    }
    \label{main_asymptotic}
\end{theorem}
\begin{proof}
    Let $\xi\in C^1(\bar\Omega)$ be such that $\nabla \xi\cdot \eta = 0$ on $\partial \Omega$, and $\rho\in C_c^1(\mathbb{R})$ be such that supp $\rho\subset [0,1]$.
    Next, for any fixed $n\in\mathbb{N}$, we let $v(x,t) = \rho(t-t_n)\xi(x)$. Furthermore, we assume that $\rho$ satisfies the following conditions         
        \begin{align}
            \rho(t-t_n) &:=(\ I^{\alpha}_{t|t_n+1} \  \hat{\rho}'\ )(t-t_n),
            \quad\text{for some $\hat{\rho}$ with }
            \hat{\rho}(0)=\hat{\rho}(1);
            \label{eq: rho}
            \\
            \tilde{\rho}(t-t_n)&:=(\ I^{1-\alpha}_{t|t_n+1} \ \rho\ )(t-t_n),
            \qquad
            \tilde{\rho}(0)=\tilde{\rho}(1)=0.
            \label{eq: rho tilde}
        \end{align}
    Multiplying equation \eqref{Eq1} by $v$ and integrating over  $(t_n-1,t_n+1)\times \Omega$ gives
        \begin{align*}
            &\mathcal{I}+\mathcal{J}=\mathcal{K},
            \intertext{where, }
            \mathcal{I} &:=\int_{t_n}^{t_n+1}\int_\Omega v\ \mathcal{D}^{\alpha}_{0|t}u(t)\ dx\ dt;
            \\
            \mathcal{J} &:= \int_{t_n}^{t_n+1}\int_\Omega v\ (-\Delta_N)^\sigma u\ dx\ dt ;
            \\
            \mathcal{K} &:= \int_{t_n}^{t_n+1}\int_\Omega v\ f(u)\ dx\ dt.
        \end{align*}
    
    \noindent First, we prove that the first term in the left-hand side, i.e $\mathcal{I}$ tends to zero as $n\to \infty$. 
    Observe that using fractional integration by parts given in Lemma \ref{lemma: int by parts}, we have 
        \begin{align}
            \mathcal{I}&=\int_{t_n}^{t_n+1} \rho(t-t_n)\ \mathcal{D}^{\alpha}_{0|t} u(t)\ dt
            \nonumber
            \\ 
            &= \int_{t_n}^{t_n+1} u(t)\ {}^{RL}D^{\alpha}_{t|t_n+1} \rho(t-t_n)\ dt\ + \left[ u(t) \, I_{t|t_n+1}^{1 - \alpha} \rho(t-t_n) \right]_{t = t_n}^{t = t_n+1}.
            \nonumber
            \intertext{From \eqref{eq: rho tilde}, we can see that $\rho$ is chosen such that $\tilde{\rho}(t-t_n)=(I_{t|t_n+1}^{1 - \alpha} \rho)(t-t_n)$ satisfies $\tilde{\rho}(t_n)=\tilde{\rho}(t_n+1)=0$, thus we have} 
            \mathcal{I}
            &= \int_{t_n}^{t_n+1} u(t)\ {}^{RL}D^{\alpha}_{t|t_n+1} \rho(t-t_n)\ dt.
            \nonumber
            \intertext{
            Notice that, 
            $({}^{RL}D^{\alpha}_{t|t_n+1} \rho)(t-t_n)
            =({}^{RL}D^{\alpha}_{t|t_n+1} I^{\alpha}_{t|t_n+1} \hat{\rho}')(t-t_n)
            =\hat{\rho}'(t-t_n)$. Thus, we have}
            \mathcal{I}
            &= \int_{t_n}^{t_n+1} u(t)\ ({}^{RL}D^{\alpha}_{t|t_n+1} I^{\alpha}_{t|t_n+1} \hat{\rho}')(t-t_n)\ dt
            \nonumber
            \\
            &=\int_{t_n}^{t_n+1} u(t)\ \hat{\rho}'(t-t_n)\ dt
            \nonumber
            \\
            &=\int_{0}^{1} u(\tau+t_n)\ \hat{\rho}'(\tau)\ d\tau,
            \label{ineq LDCT asymp proof term 1}
        \end{align}
      where we have used the change of variable $\tau = t-t_n$. Since $u(x,t)$ is uniformly continuous for $t\in [\delta, \infty)$ for $\delta>0$, by Lebesgue's dominated convergence theorem, we conclude that the integrand of the right-hand side of \eqref{ineq LDCT asymp proof term 1} is bounded and goes to zero for large $n$. Indeed, 
            \begin{align*}
                \lim_{n\to\infty} \int_{0}^{1} u(\tau+t_n)\ \hat{\rho}'(\tau)\ d\tau = u_\infty (\hat{\rho}(1)-\hat{\rho}(0))=0.
            \end{align*}
    
    For the remaining part of the equation, using Lemma \ref{Lemma_int by parts lap},  we notice that $\mathcal{J}$ can be written as
        \begin{align}
            \mathcal{J}=\int_{t_n}^{t_n+1}\rho(t-t_n)\int_\Omega u\ (-\Delta_N )^\sigma\xi(x)\ dx\ dt.
        \end{align}
    Hence, $\mathcal{J}+\mathcal{K}=0$ is equivalent to
        \begin{align*} 
            \int_{-1}^{1}\rho(\tau)\ d\tau \bigg(\int_\Omega u(x,\tau+t_n)\ (-\Delta_N)^\sigma \xi(x)
            &- \xi(x)f(u(x,\tau+t_n))\ dx \bigg) 
            =0.
            \intertext{Using Lebesgue's theorem, taking the limit as $n\to\infty$, we have} 
            \int_{-1}^{1}\rho(\tau) \ d\tau \bigg(\int_\Omega u_{\infty}\ (-\Delta_N )^\sigma \xi(x)&- \xi(x) f(u_\infty)\ dx \bigg)\
            =0,
        \end{align*}
    where we have used the uniform continuity result of Theorem \ref{uniform_cont thm}. Consequently, using \eqref{eq: rho}, we get
        \begin{align*}
            \int_\Omega u_{\infty}\ (-\Delta_N )^\sigma \xi(x) &- \xi(x)f(u_\infty) \ dx \
            =0,
            \intertext{or}
            \int_\Omega \xi(x)\ (-\Delta_N )^\sigma u_{\infty}&- \xi(x)f(u_\infty) \ dx \
            =0.
        \end{align*}
    which is the variational formulation of the elliptic problem \eqref{eq: stationary}.

    To show that $u$ tends to $u_\infty= 1/\sqrt{b}$, we multiply \eqref{eq: stationary} by $(1-\sqrt{b}u_\infty)$ and integrate over $\Omega$ to get
        \begin{align*}
           \int_{\Omega} (-\Delta_N)^\sigma u_\infty \ (1-\sqrt{b}u_\infty)\ dx
           &=\int_{\Omega} u_\infty(1+\sqrt{b}u_\infty)(1-\sqrt{b}u_\infty)^2\ dx.\\
           \intertext{Using integration by parts, we get}
           \int_{\Omega} (-\Delta_N)^{\sigma/2} u_\infty \ (-\Delta_N)^{\sigma/2}  (1-\sqrt{b}u_\infty) \ dx
           &=\int_{\Omega} u_\infty(1+\sqrt{b}u_\infty)(1-\sqrt{b}u_\infty)^2\ dx,
           \intertext{or,}
            -\sqrt{b}\int_{\Omega} |(-\Delta_N)^{\sigma/2} u_\infty|^2 \ dx
           &=\int_{\Omega} u_\infty(1+\sqrt{b}u_\infty)(1-\sqrt{b}u_\infty)^2\ dx.
        \end{align*}
    As we have established the positivity of our solution, we must have
        \begin{align*}
            \int_{\Omega} |(-\Delta_N)^{\sigma/2} u_\infty|^2 \ dx=0
            \quad\text{and}\quad
            \int_{\Omega} u_\infty(1+\sqrt{b}u_\infty)(1-\sqrt{b}u_\infty)^2\ dx=0.
        \end{align*}
    This means
        \begin{align*}
            (-\Delta_N)^{\sigma/2} u_\infty=0
            \qquad\Longrightarrow\qquad
            u_\infty \text{ is a constant},
        \end{align*}
    and from $\int_{\Omega} u_\infty(1+\sqrt{b}u_\infty)(1-\sqrt{b}u_\infty)^2\ dx=0$, we have
         \begin{align*}
            u_\infty(1+\sqrt{b}u_\infty)(1-\sqrt{b}u_\infty)^2=0
            \quad\Longrightarrow\quad
            u_\infty=0 
            \quad\text{or}\quad
            u_\infty=1/\sqrt{b}.
        \end{align*}
    Lastly, we show that $u_\infty=0$ is unstable. Let $u=u_\infty+\bar{u}=\bar{u}$. Substituting in \eqref{Eq1} we get the following linearized equation
        \begin{align}
             \mathcal{D}_{a|t}^{\alpha} \bar{u} + (-\Delta_N)^{\sigma/2}\bar{u} = \bar{u}. 
             \label{eq: linearised}
        \end{align}
    Whose solution is given by 
        \begin{align*}
            \bar{u}(x,t) = E_\alpha\left(t^\alpha\mathcal{L}\right)u_0(x);
            \qquad \mathcal{L}= I-(-\Delta_N)^{\sigma/2}.
        \end{align*}
    If there exist $\lambda\in spec(\mathcal{L})$ with $\arg(\lambda)\le\alpha\pi/2$ then $\bar{u}(x,t) \not\to 0$  or $u_\infty=0$ is unstable. In our case, it is sufficient to show the existence of positive $\lambda\in spec(\mathcal{L})$.
    
    To do that consider the eigenvalue problem
        \begin{align}
            \left(1-(-\Delta_N)^{\sigma/2}\right) \psi = \lambda\psi.
            \label{eq: eigenvalue problem}
        \end{align}
    If there exist $\lambda>0$ such that $(1-\lambda)>0$ then the zero solution is unstable. Indeed it is the case; 
    we see that $\lambda=1$ is the eigenvalue that corresponds to the eigenfunction $\psi=const$ of the problem \eqref{eq: eigenvalue problem}.
    Hence our assertion follows.    
\end{proof}

\section*{Acknowledgments}

The authors would like to express their gratitude to Khalifa University of Science and Technology for their continuous support and resources that made this research possible.

\end{document}